\documentclass{amsart}
\usepackage{graphicx}

\newtheorem{theorem}{Theorem}[section]

\newtheorem{corollary}[theorem]{Corollary}
\newtheorem{proposition}[theorem]{Proposition}
\newtheorem{observation}[theorem]{Observation}

\theoremstyle{definition}

\newtheorem{example}[theorem]{Example}
\newtheorem{xca}[theorem]{Exercise}
\newtheorem{claim}[theorem]{Claim}

\theoremstyle{remark}
\newtheorem{remark}[theorem]{Remark}

\numberwithin{equation}{section}

\begin{document}

\title{A property of diagrams of the trivial knot}

\author{Makoto Ozawa}
\address{Department of Natural Sciences, Faculty of Arts and Sciences, Komazawa University, 1-23-1 Komazawa, Setagaya-ku, Tokyo, 154-8525, Japan}
\email{w3c@komazawa-u.ac.jp}





\begin{abstract}
In this paper, we give a necessary condition for a diagram to represent the trivial knot.
\end{abstract}

\maketitle

\section{Introduction}

\subsection{What can we say when a diagram represents the trivial knot?}


Let $K$ be a knot in the 3-sphere $S^3$ and consider a diagram $\pi(K)$ of $K$ on the 2-sphere $S^2$.
We say that a diagram $\pi(K)$ is {\it $\rm{I}$-reduced} (resp. {\it $\rm{II}$-reduced}) if the crossing number of $\pi(K)$ cannot be reduced by a Reidemeister move $\rm{I}$ (resp. Reidemeister move $\rm{II}$).
We say that a diagram $\pi(K)$ is {\it prime} if it contains at least one crossing and for any loop $l$ intersecting $D$ in two points except for crossings, there exists a disk $D$ in $S^2$ such that $D\cap \pi(K)$ consists of an embedded arc.

\begin{figure}[htbp]
	\begin{center}
		\includegraphics[trim=0mm 0mm 0mm 0mm, width=.5\linewidth]{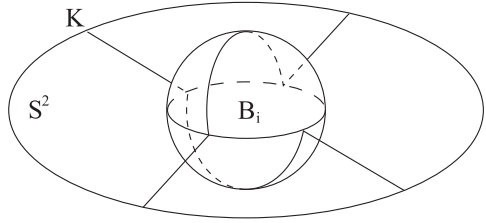}
	\end{center}
	\caption{a bubble between the over arc and the under arc}
	\label{bubble1}
\end{figure}

\begin{figure}[htbp]
	\begin{center}
		\includegraphics[trim=0mm 0mm 0mm 0mm, width=.5\linewidth]{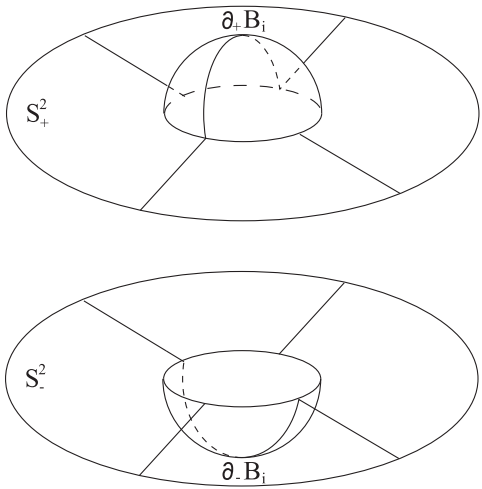}
	\end{center}
	\caption{the upper hemisphere and the lower hemisphere}
	\label{bubble2}
\end{figure}

We position $K$ in the following Menasco's manner (\cite{M}).
For each crossing $c_i$ of $\pi(K)$, we insert a small 3-ball ``bubble" $B_i$ between the over crossing and the under crossing of $c_i$ and isotope the over arc of $N(c_i;K)$ onto the upper hemisphere $\partial_+ B_i$ of $\partial B_i$ and the under arc onto the lower hemisphere $\partial_- B_i$.
See Figure \ref{bubble1}.
Let $S^2_+$ (resp. $S^2_-$) be a 2-sphere obtained from $S^2$ by replacing each equatorial disk $B_i\cap S^2$ with the upper (resp. lower) hemisphere $\partial_+ B_i$ (resp. $\partial_- B_i$).
See Figure \ref{bubble2}.
Put $P=S^2_+\cap S^2_-$.
Then, $K$ is contained in $S^2_+\cup S^2_-=P\cup \bigcup \partial B_i$.
We call each component of $P-\pi(K)$ a {\it region}.
Let $B$ be the union of all bubbles and $R$ be the union of all regions.

A loop $l$ embedded in $S^2_+-K$ (resp. $S^2_--K$) is called a {\it $+$-Menasco loop} (resp. {\it $-$-Menasco loop}) if for each region $R_j$, each component of $l\cap R_j$ is an arc connecting different arc components of $\partial B\cap R_j$ and for each bubble $B_i$, each component of $l\cap \partial B_i$ is an arc connecting two different regions.
The number of crossings which a Menasco loop passes through is called {\it length}.

\begin{xca}
Show that there is always a length two $+$ or $-$-Menasco loop if a diagram is not prime and contains at least one crossing.
\end{xca}

\begin{figure}[htbp]
	\begin{center}
		\includegraphics[trim=0mm 0mm 0mm 0mm, width=.6\linewidth]{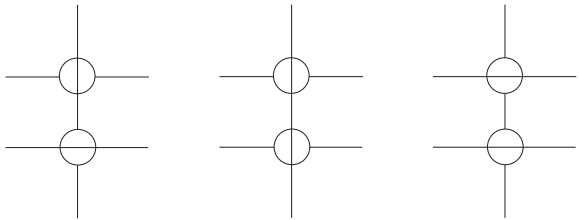}
	\end{center}
	\caption{adjacent, $+$-adjacent, $-$-adjacent}
	\label{adjacent}
\end{figure}

Two crossings $c_i$ and $c_j$ are {\it adjacent} if there exists an arc $\gamma$ of $K\cap P$ connecting two bubbles $B_i$ and $B_j$, and {\it $+$-adjacent} (resp. {\it $-$-adjacent}) if $\gamma$ connects two over arcs $K\cap \partial_+ B_i$ and $K\cap \partial_+ B_j$ (resp. two under arcs $K\cap \partial_- B_i$ and $K\cap \partial_- B_j$).
See Figure \ref{adjacent}.

\begin{xca}
Show that if there exists a length two $+$-Menasco loop passing through two $+$-adjacent crossings, then the diagram is not prime or not $\rm{II}$-reduced.
\end{xca}



The following is a restriction on Corollary \ref{trivial knot} which is an essence of this paper.

\begin{corollary}\label{trivial}
Any $\rm{I}$-reduced, $\rm{II}$-reduced, prime diagram of the trivial knot has a $\pm$-Menasco loop passing through $2n$-crossings $c_1,c_2,\ldots,c_{2n}$, where $n\ge 2$ and $c_{2i-1}$ is $\pm$-adjacent to $c_{2i}$ for $i=1,\ldots,n-1$.
\end{corollary}

\begin{example}
Consider a $\rm{I}$-reduced, $\rm{II}$-reduced, prime, 4-crossing diagram of the right-handed trefoil.
See Figure \ref{example1}.
There exists a $+$-Menasco loop passing through $c_1, c_2, c_4, c_3, c_1, c_3$, where $c_1$ is $+$-adjacent to $c_2$ and $c_4$ is $+$-adjacent to $c_3$.
After Theorem \ref{main}, we will know that any $\rm{I}$-reduced, $\rm{II}$-reduced diagram of the trefoil except for the 3-crossing diagram has a $\pm$-Menasco loop satisfying the condition in Corollary \ref{trivial}.

\begin{figure}[htbp]

	\begin{center}
		\includegraphics[trim=0mm 0mm 0mm 0mm, width=.3\linewidth]{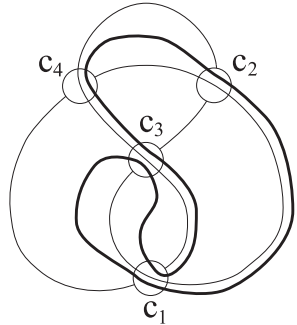}
	\end{center}
	\caption{$\rm{I}$-reduced, $\rm{II}$-reduced, prime, 4-crossing diagram of the right-handed trefoil}
	\label{example1}
\end{figure}

\end{example}

\begin{example}
Next example is borrowed from Ochiai's book \cite{O2}.
This diagram of the trivial knot has no $r$-wave for any $r\ge 0$.
See Figure \ref{example2}.
At each stage, there exists a $\pm$-Menasco loop satisfying the condition in Corollary \ref{trivial} or it is not $\rm{I}$-reduced or not $\rm{II}$-reduced.
In the former case, a $\pm$-Menasco loop can be used to simplify the diagram if it has successive three adjacent crossings, and in the latter case, the crossing number can be reduced by a Reidemeister move $\rm{I}$ or $\rm{II}$.

\begin{figure}[htbp]
	\begin{center}
		\includegraphics[trim=0mm 0mm 0mm 0mm, width=.8\linewidth]{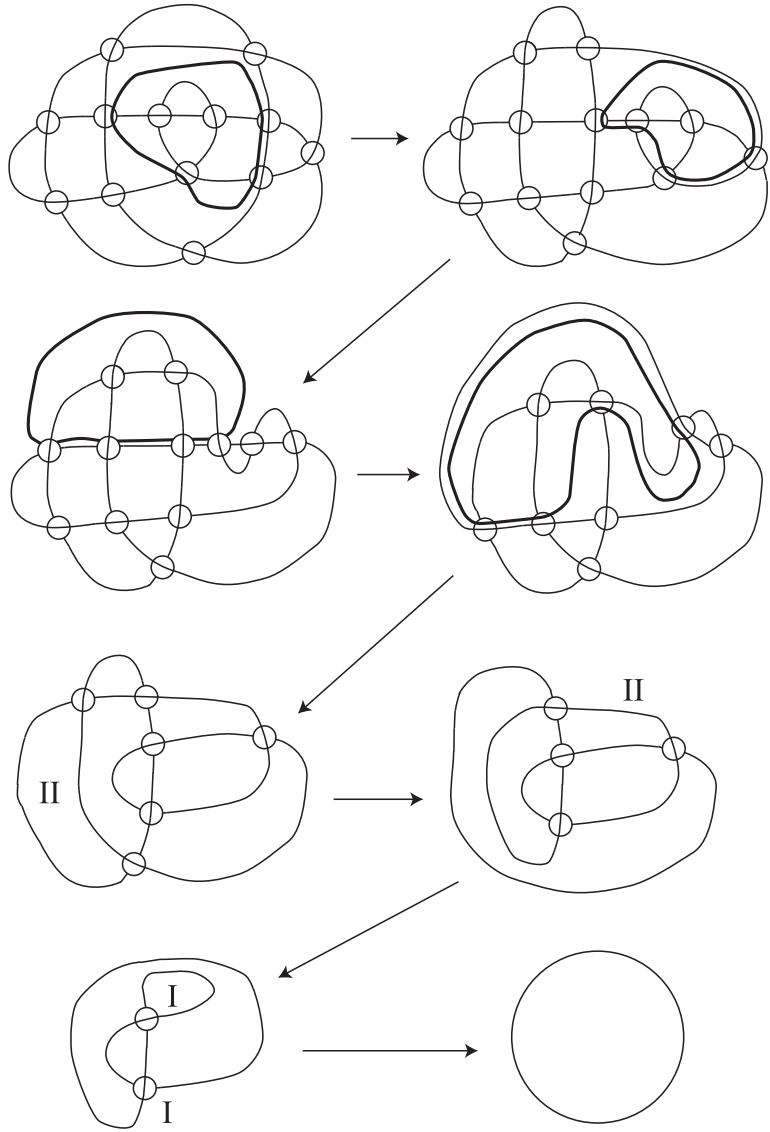}
	\end{center}
	\caption{Ochiai's ``non-trivial'' diagram of the trivial knot}
	\label{example2}
\end{figure}


\end{example}

\begin{example}
Final example is somewhat artificial.
This diagram is $2$-almost alternating, that is, obtained from an alternating diagram by twice crossing changes on it.
See Figure \ref{example3}.
There does not exist a $\pm$-Menasco loop satisfying the condition in Corollary \ref{trivial}.
Hence, this knot is non-trivial.

Note that Tsukamoto characterized almost alternating diagarms of the trivial knot (\cite{T}).
\begin{figure}[htbp]
	\begin{center}
		\includegraphics[trim=0mm 0mm 0mm 0mm, width=.6\linewidth]{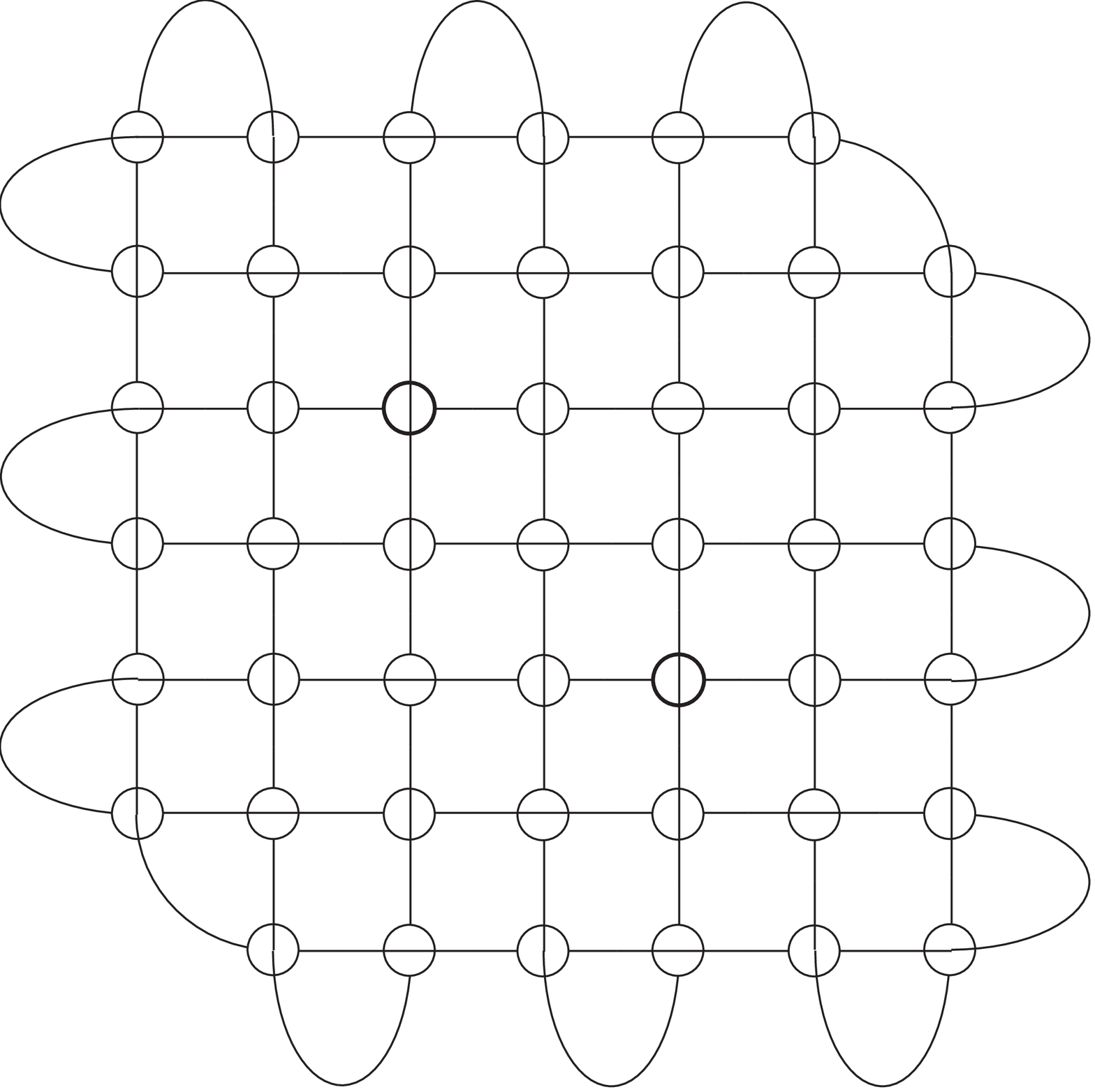}
	\end{center}
	\caption{2-almost alternating diagram}
	\label{example3}
\end{figure}

\end{example}

\subsection{Where do you untie it from?}


We position a knot $K$ in the Menasco's manner mentioned above.
In \cite{MT}, Menasco and Thistlethwaite defined a standard position of a spanning surface $F$ for $K$ as follows.

\begin{enumerate}
	\item[(i)] $\rm{int}F$ meets each of $S_+,S_-$ transversely in a pairwise disjoint collection of simple closed curves and arcs;
	\item[(ii)] $F$ meests each $B_i$ in a collection of saddle-shaped disks;
	\item[(iii)] there is a collar $C\cong I\times \partial F$ of $\partial F$ in $F$ and a projection $p:C\to \partial F$ such that for each $x\in \partial F\cap \partial B_i$ the fibre $p^{-1}(x)$ is a straight line segment which is normal to $\partial B_i$ and which does not meet the interior of $B_i$.
\end{enumerate}

Moreover they showed;

\begin{proposition}[Proposition 2 in \cite{MT}]
The disk $F$ spanning $K$ may be replaced with a spanning disk $F'$ such that each circle $C$ in $F'\cap S_+$ satisfies the following:
\begin{enumerate}
	\item[(i)] $C$ bounds a disk in $F'$ whose interior lies entirely above $S_+$;
	\item[(ii)] if $C\subset \rm{int}F'$, $C$ meets at least one bubble;
	\item[(iii)] $C$ does not meet any bubble in more than one arc $($whether or not $C\subset \rm{int}F'$$)$.
\end{enumerate}
Moreover, the corresponding conditions for $F\cap S_-$ can be achieved simultaneously.
\end{proposition}

\begin{figure}[htbp]
	\begin{center}
		\includegraphics[trim=0mm 0mm 0mm 0mm, width=.7\linewidth]{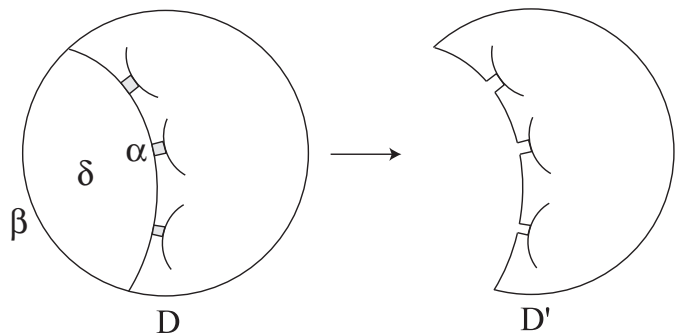}
	\end{center}
	\caption{an outermost Menasco arc}
	\label{outermost arc}
\end{figure}

Hereafter, let $K$ be the trivial knot and $D$ be a disk bounded by $K$.
We assume the above conditions on $D$.
Then, there exists an outermost arc $\alpha$ on $D$ which bounds an outermost disk $\delta$ in $D$ with a subarc $\beta$ of $K$.
We call an arc of $D\cap S_{\pm}$ a {\it $\pm$-Menasco arc}.
See Figure \ref{outermost arc}.

\begin{figure}[htbp]
	\begin{center}
		\includegraphics[trim=0mm 0mm 0mm 0mm, width=.7\linewidth]{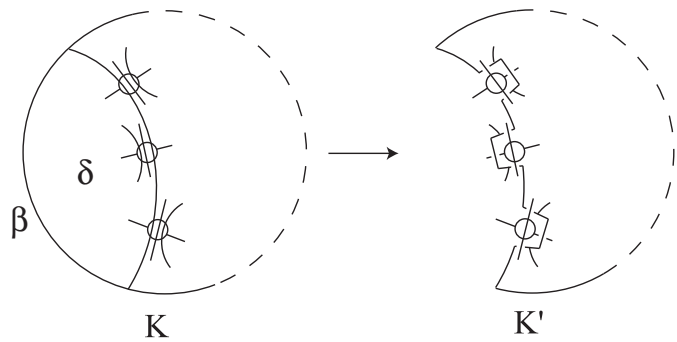}
	\end{center}
	\caption{Moving $K$ along the outermost Menasco arc}
	\label{outermost arc2}
\end{figure}

Let $\delta'$ be a subdisk of $D$ which forms $\delta$ together with all saddle-shaped disks meeting $\delta$.
We say that $K'$ is obtained by a {\it move along an outermost $\pm$-Menasco arc} if $K'$ is obtained by isotoping $K$ along $\delta'$.
See Figure \ref{outermost arc2}.
Note that $K'$ is not in the Menasco's position, but there exists a next outermost $\pm$-Menasco arc $\alpha'$ on $D'$ and we can move $K'$ along $\alpha'$.

We can summarize our observation as follows.

\begin{observation}
Let $K_1$ be the trivial knot and $\pi(K_1)$ be a diagram of $K_1$ with at least one crossing.
Then, there exists a sequence $\alpha_1,\ldots,\alpha_n$ of $\pm$-Menasco arc such that $K_{i+1}$ is obtained from $K_i$ by a move along $\alpha_i$ and $K_{n+1}$ is a diagram without crossings.
\end{observation}

\section{Main Theorem}

In this secion, we consider diagrams of a knot on a closed surface.

Let $F$ be a closed surface embedded in $S^3$ and $K$ a knot contained in $F\times [-1,1]$.
Suppose that $\pi(K)$ is a regular projection on $F$, where $\pi : F\times [-1,1]\to F\times \{0\}=F$ is the projection.
Then, we have a regular diagram on $F$ obtained from $\pi(K)$ by adding the over/under information to each double point, and we denote it by the same symbol $\pi(K)$ in this article.

We say that a diagram $\pi(K)$ on $F$ is {\em reduced} if there is no disk region of $F-\pi(K)$ which meets only one crossing.
We say that a diagram $\pi(K)$ on $F$ is {\em prime} if it contains at least one crossing and for any loop $l$ intersecting $\pi(K)$ in two points except for crossings, there exists a disk $D$ in $F$ such that $D\cap \pi(K)$ consists of an embedded arc.

Let $S$ be a closed surface of positive genus in $S^3$ and $K$ a knot contained in $S$.
The {\em representativity} $r(S,K)$ of a pair $(S,K)$ is defined as the minimal number of intersecting points of $K$ and $\partial D$, where $D$ ranges over all compressing disks for $S$ in $S^3$.
It follows from Lemma 3 in \cite{MO} that $r(S,K)\ge1$ if and only if $S\cap E(K)$ is incompressible in $E(K)$, and $r(S,K)\ge2$ if and only if $S\cap E(K)$ is incompressible and $\partial$-incompressible in $E(K)$, where $E(K)$ denotes the exterior of $K$ in $S^3$.

In the previous paper, the author showed the non-triviality of generalized alternating knots by the following theorem.

\begin{theorem}[\cite{MO}]\label{alternating}
Let $F$ be a closed surface embedded in $S^3$, $K$ a knot contained in $F\times [-1,1]$ which has a reduced, prime, alternating diagram on $F$.
Then, we have the following.
\begin{enumerate}
	\item $F-\pi(K)$ consists of open disks.
	\item $F-\pi(K)$ admits a checkerboard coloring.
	\item $K$ bounds a non-orientable surface $H$ coming from the checkerboard coloring.
	\item $K$ can be isotoped into $\partial N(H)$ so that $\partial N(H)-K$ is connected.
	\item $r(\partial N(H),K)\ge 2$.
\end{enumerate}
\end{theorem}

We call the closed surface $\partial N(H)$ an {\it interpolating surface} obtained from the checkerboard coloring, where $H$ is one of the checkerboard surfaces.

Theorem \ref{alternating} assures the existence of an incompressible and $\partial$-incompressible separating orientable surface of integral boundary slope in the exterior of a generalized alternating knot.
Hence, the knot is non-trivial.

Conversely, suppose that an interpolating surface obtained from the checkerboard coloring is compressible.
Then, we know from Theorem \ref{alternating} that the diagram is not alternating.
Are there other properties of the diagram?
This is the main subject in this paper.

\begin{figure}[htbp]
	\begin{center}
		\includegraphics[trim=0mm 0mm 0mm 0mm, width=.7\linewidth]{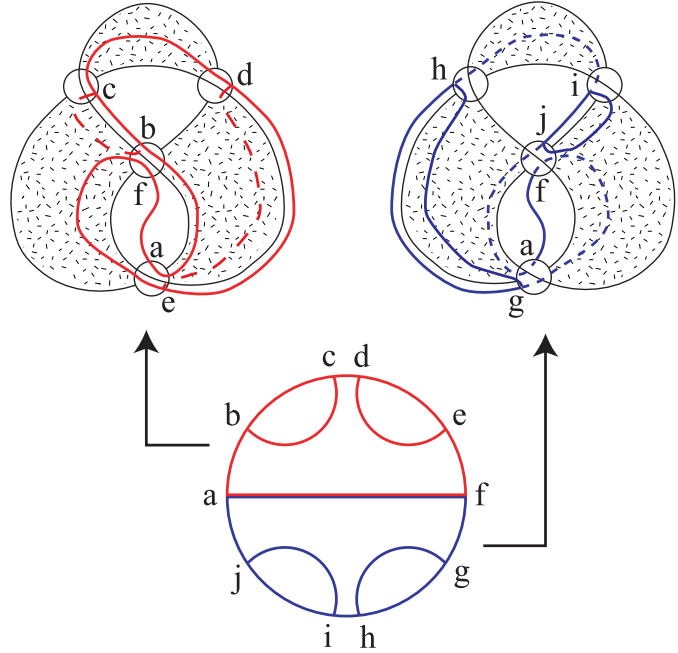}
	\end{center}
	\caption{a compressing disk for the interpolating surface}
	\label{comp_check}
\end{figure}

\begin{example}
The 4-crossing diagram of the right-handed trefoil without nugatory crossings.
The interpolating surface obtained from a checkerboard surface is compressible.
See Figure \ref{comp_check}.
Here, the compressing disk intersects the union of regions in 5 arcs.
\end{example}

As Introduction, we position $K$ in the Menasco's manner with respect to a closed surface $F$.
For each crossing $c_i$ of $\pi(K)$, we insert a small 3-ball ``bubble" $B_i$ between the over crossing and the under crossing of $c_i$ and isotope the over arc of $N(c_i;K)$ onto the upper hemisphere $\partial_+ B_i$ of $\partial B_i$ and the under arc onto the lower hemisphere $\partial_- B_i$.
Let $F_+$ (resp. $F_-$) be a 2-sphere obtained from $F$ by replacing each equatorial disk $B_i\cap F$ with the upper (resp. lower) hemisphere $\partial_+ B_i$ (resp. $\partial_- B_i$).
Put $P=F_+\cap F_-$.
Then, $K$ is contained in $F_+\cup F_-=P\cup \bigcup \partial B_i$.
We call each component of $P-\pi(K)$ a {\it region}.
Let $B$ be the union of all bubbles and $R$ be the union of all regions.

The following is a main theorem.

\begin{theorem}\label{main}
Let $F$ be a closed surface embedded in $S^3$, $K$ a knot contained in $F\times [-1,1]$ which has a $\rm{I}$-reduced, $\rm{II}$-reduced, prime, checkerboard colorable diagram on $F$.
If at least one of two interpolating surfaces obtained from the checkerboard coloring is compressible in the complement of $K$, then there exists a compressing disk $\delta$ for $F_{\pm}-K$ in $V_{\pm}$ such that $\partial \delta$ is a $\pm$-Menasco loop passing through $2n$-crossings $c_1,c_2,\ldots,c_{2n}$, where $n\ge 2$ and $c_{2i-1}$ is $\pm$-adjacent to $c_{2i}$ for $i=1,\ldots,n-1$.
\end{theorem}

\begin{remark}
Theorem \ref{main} still holds for link case.
If a link is split but the diagram is not connected, then both of checkerboard surfaces are compressible in the link complement and there exists a $\pm$-Menasco loop satisfying the condition in Theorem \ref{main}.
\end{remark}

\begin{remark}\label{once}
In Theorem \ref{main}, we can take a compressing disk $D$ so that $\partial D$ does not pass through a one side of a crossing {\em more than once}.
\end{remark}

\begin{xca}
Show the statement of Remark \ref{once}.
\end{xca}

\begin{remark}
It is possible to state that for a checkerboard surface $F$, whether $\tilde{F}$ is compressible by means of {\em all $\pm$-Menasco loop} coming from {\em all subdisk} in $D$.
\end{remark}

\section{Proof}

\begin{proof} (of Theorem \ref{main})
The following claim is Claim 6 in \cite{MO}.

\begin{claim}\label{open disk}
$F-\pi(K)$ consists of open disks.
\end{claim}


Let $\partial N(H)$ be an interpolating surface obtained from the ckeckerboard coloring such that $\partial N(H)-K$ is compressible in $S^3-K$.
The following claim is Claim 9 in \cite{MO}.

\begin{claim}\label{I-bundle}
$\partial N(H)-K$ is incompressible in $N(H)$.
\end{claim}

Hence, $\partial N(H)-K$ is compressible in the outside of $N(H)$, and let $D$ be a compressing disk for $\partial N(H)-K$.

We regard $N(H)$ as the following.
For each crossing $c_i$ of $\pi(K)$, we insert a small 3-ball $B_i$ as a regular neighborhood of $c_i$.
In the rest of these 3-balls, we consider the product $R_i\times I$ for each region $R_i$ of $F-\pi(K)$.
Then, the union of $B_i$'s and $R_i\times I$'s is homeomorphic to $N(N)$.
See Figure \ref{N(H)}.

\begin{figure}[htbp]
	\begin{center}
		\includegraphics[trim=0mm 0mm 0mm 0mm, width=.6\linewidth]{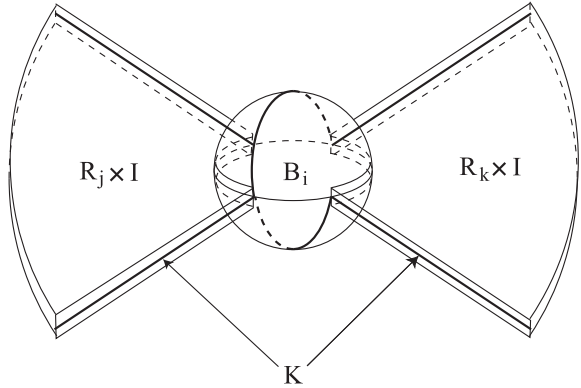}
	\end{center}
	\caption{$B_i$, $R_i\times I$ and $K$}
	\label{N(H)}
\end{figure}

Let $\Delta=\Delta_1\cup \cdots \cup \Delta_n$ be the union of components of $F-int N(H)$, where $\Delta_i$ is a disk by Claim \ref{open disk}.
Then, each component of $(\partial N(H)-K)-\partial \Delta$ is an open disk containing $R_i\times\{0\}$ or $R_i\times\{1\}$ for some $i$ and whose closure is denoted by $R_i^-$ or $R_i^+$ respectively.
Put $R=(\bigcup_i R_i^-)\cup (\bigcup_i R_i^+)$ and $B=\bigcup_i B_i$.
Note that $D\cap \Delta\ne \emptyset$, otherwise $\partial D$ is entirely contained in a region and $D$ would not be a compressing disk for $\partial N(H)-K$.

The following claim is Claim 10 in \cite{MO}.

\begin{claim}\label{standard}
We may assume the following.
\begin{enumerate}
	\item $\partial D\cap R$ consists of arcs that connect different arc components of $\partial B\cap \partial \Delta$.
	\item $D\cap \Delta$ consists of arcs that connect different arc components of $\partial B\cap \partial \Delta$.
\end{enumerate}
\end{claim}

Next, we concentrate on an outermost arc $\alpha$ of $D\cap \Delta$ in $D$ and the corresponding outermost disk $\delta$ in $D$.
Put $\delta\cap \partial D=\beta$.

\begin{claim}\label{outermost}
Any outermost arc $\alpha$ connects $\pm$-adjacent crossings.
\end{claim}

\begin{proof}

By Claim \ref{standard}, we have two configurations.

\begin{description}
	\item[Case 1] $\beta$ connects the same crossing ball $B_i$ (Figure \ref{alternate1}).
	\item[Case 2] $\beta$ connects different crossing balls $B_i$ and $B_j$ (Figure \ref{alternate2}).
\end{description}

\begin{figure}[htbp]
	\begin{center}
		\includegraphics[trim=0mm 0mm 0mm 0mm, width=.6\linewidth]{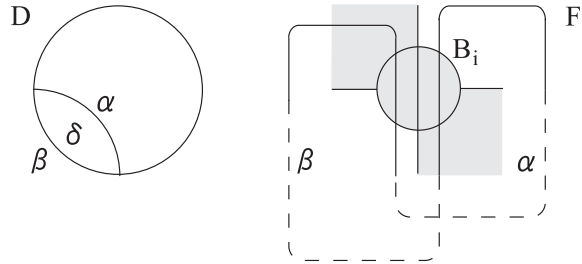}
	\end{center}
	\caption{Configuration of Case 1}
	\label{alternate1}
\end{figure}

\begin{figure}[htbp]
	\begin{center}
		\includegraphics[trim=0mm 0mm 0mm 0mm, width=.9\linewidth]{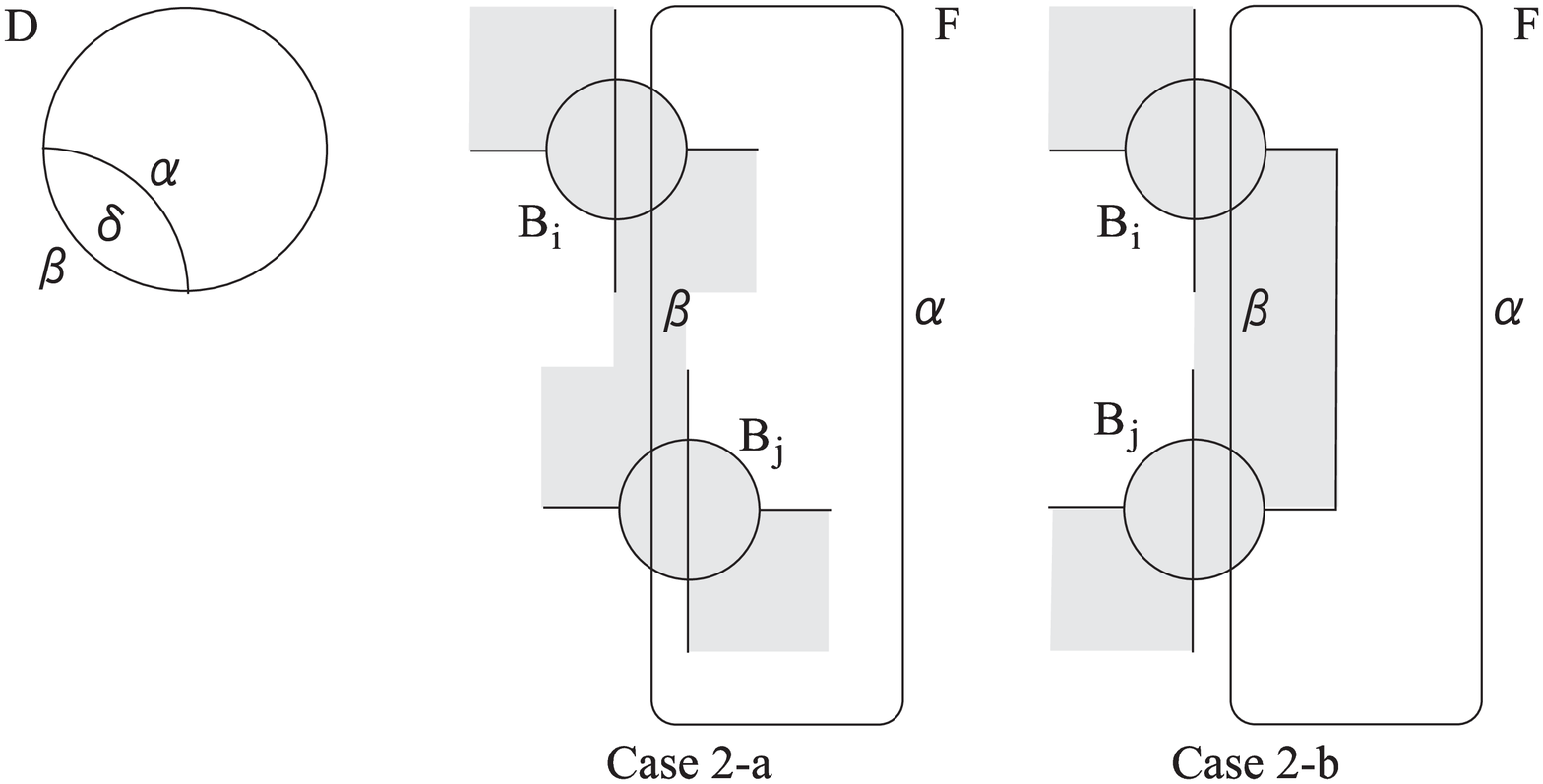}
	\end{center}
	\caption{Configuration of Case 2}
	\label{alternate2}
\end{figure}

In Case 1, by connecting $\partial \beta$ on $\partial B_i$ and projecting on $F$, we have a loop $l_\beta$ on $F$ which intersects $\pi(K)$ in one crossing point $c_i$.
Similarly, we obtain a loop $l_\alpha$ on $F$ which intersects $\pi(K)$ in one crossing point $c_i$.
Since $l_\beta$ intersects $l_\alpha$ in one point $c_i$, $l_\beta$ is essential in $F$.
Let $l_\beta$ avoid $c_i$.
Then we have an essential loop in $F$ which intersects $\pi(K)$ in two points except for crossings.
This contradicts that $\pi(K)$ is prime.

In Case 2, we have a loop $\pi(\alpha\cup\beta)$ in $F$ which intersects $\pi(K)$ in two points except for crossings.
In Case 2-a, it does not bound a disk $D'$ in $F$ such that $D'\cap \pi(K)$ is an arc since there are crossings $c_i$ and $c_j$ on both sides of the loop.
This contradicts that $\pi(K)$ is prime.
In Case 2-b, the loop $\pi(\alpha\cup\beta)$ bounds a disk $D'$ such that $D'\cap \pi (K)$ is an arc as in Figure \ref{alternate2} since $\pi(K)$ is prime.
This shows that $B_i$ and $B_j$ are $\pm$-adjacent.
\end{proof}

If $|D\cap \Delta|=1$, then $\pi(K)$ is $\rm{II}$-reducible.
See Figure \ref{II-reducible}.
This contradicts the supposition of Theorem \ref{main}.

\begin{figure}[htbp]
	\begin{center}
		\includegraphics[trim=0mm 0mm 0mm 0mm, width=.5\linewidth]{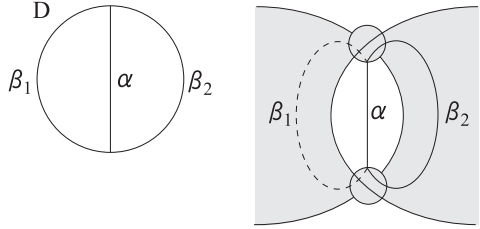}
	\end{center}
	\caption{$\rm{II}$-reducible}
	\label{II-reducible}
\end{figure}

Otherwise, there exists an outermost fork in the graph on $D$ obtained by $D\cap \Delta$, where we regard the closure of each component of $D-\Delta$ as a vertex, and connects vertices if the corresponding components are close to each other.
See Figure \ref{outermost fork}.
Then, a region $\delta$ corresponding to an outermost fork gives a compressing disk for $F_{\pm}-K$ in $V_{\pm}$ such that $\partial \delta$ is a $\pm$-Menasco loop passing through $2n$-crossings $c_1,c_2,\ldots,c_{2n}$.
Claim \ref{outermost} satisfies the condition that $n\ge 2$ and $c_{2i-1}$ is $\pm$-adjacent to $c_{2i}$ for $i=1,\ldots,n-1$.
See Figure \ref{M-loop}.
\begin{figure}[htbp]
	\begin{center}
		\includegraphics[trim=0mm 0mm 0mm 0mm, width=.35\linewidth]{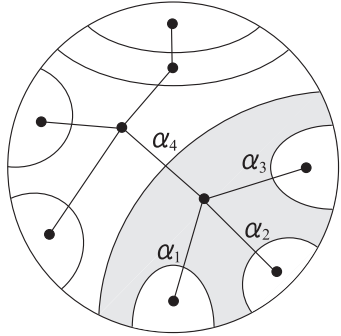}
	\end{center}
	\caption{outermost fork}
	\label{outermost fork}
\end{figure}
\begin{figure}[htbp]
	\begin{center}
		\includegraphics[trim=0mm 0mm 0mm 0mm, width=.5\linewidth]{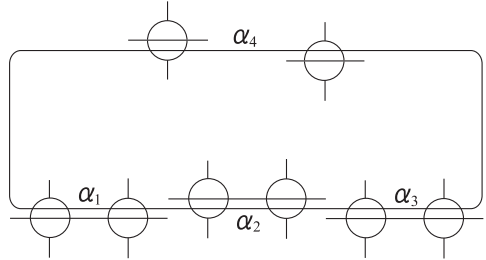}
	\end{center}
	\caption{$+$-Menasco loop}
	\label{M-loop}
\end{figure}
\end{proof}

Here, we state Corollary \ref{trivial} in more general form.

\begin{corollary}\label{trivial knot}
Let $F$ be a closed surface embedded in $S^3$, $K$ a knot contained in $F\times [-1,1]$ which has a $\rm{I}$-reduced, $\rm{II}$-reduced, prime, checkerboard colorable diagram on $F$.
If $K$ is trivial, then there exists a compressing disk $\delta$ for $F_{\pm}-K$ in $V_{\pm}$ such that $\partial \delta$ is a $\pm$-Menasco loop passing through $2n$-crossings $c_1,c_2,\ldots,c_{2n}$, where $n\ge 2$ and $c_{2i-1}$ is $\pm$-adjacent to $c_{2i}$ for $i=1,\ldots,n-1$.
\end{corollary}

\begin{proof} (of Corollary \ref{trivial knot})
Let $\pi(K)$ be a $\rm{I}$-reduced, $\rm{II}$-reduced, prime, checkerboard colorable diagram on $F$ and $S$ an interpolating surface obtained from a checkerboard surface $H$.

If $S\cap E(K)$ is incompressible and $\partial$-incompressible in $E(K)$, then by Lemma 1 in \cite{MO}, each component of $S\cap E(K)$ is a disk.
Hence $S$ is a 2-sphere and $H$ is a disk.
Then, $\pi(K)$ is $\rm{I}$-reducible or it has no crossing.
The former contradicts that $\pi(K)$ is $\rm{I}$-reduced and the latter contradicts that $\pi(K)$ is prime.

If $S\cap E(K)$ is compressible in $E(K)$, then by Theorem \ref{main}, the conclusion of Corollary \ref{trivial knot} is satisfied.

Otherwise, $S\cap E(K)$ is incompressible and $\partial$-compressible in $E(K)$.
By Lemma 2 in \cite{MO}, each component of $S\cap E(K)$ is $\partial$-parallel annulus.
Hence $S$ is a torus and $H$ is a M\"{o}bius band.
Since $\pi(K)$ is $\rm{I}$-reduced and $\rm{II}$-reduced, $\pi(K)$ is a standard $(2,n)$-torus knot diagram, where $n$ is an odd integer.
If $|n|\ge3$, then $r(S,K)=2$ and $S\cap E(K)$ is $\partial$-incompressible in $E(K)$.
This contradicts the assumption.
Otherwise, $n=\pm1$.
This shows that $\pi(K)$ is $\rm{I}$-reducible, a contradiction.
\end{proof}

\bibliographystyle{amsplain}

\end{document}